\documentclass[a4paper,12pt]{article}
\usepackage[top=2.5cm,bottom=2.5cm,left=2.5cm,right=2.5cm]{geometry}
\usepackage{cite, amsmath, amssymb}
\usepackage[margin=1cm,%
            font=small,%
            format=hang,%
            labelsep=period,%
            labelfont=bf]{caption}
\usepackage{graphicx}
\usepackage{cite}
\allowdisplaybreaks[1]

\begin{document}
\title{\vspace{3.5 cm} Ordering Trees by Their ABC Spectral Radii
\footnote{Supported by the NSFC (11771362) and the Anhui Provincial Natural Science Foundation (2008085J01).}
}
\author{Wenshui Lin$^{1,2}$, Zhangyong Yan$^2$, Peifang Fu$^2$, Jia-Bao Liu $^{3,}$\footnote{Corresponding author. E-mail address: liujiabaoad@163.com (J. B. Liu).}\\
        {\small $^1$ Fujian Key Laboratory of Sensing and Computing for Smart City, Xiamen 361005, China}\\
        {\small $^2$ School of Informatics, Xiamen University, Xiamen 361005, China}\\
        {\small $^3$ School of Mathematics and Physics, Anhui Jianzhu University, Hefei 230601, China}
}

\date{\small (Received August 03, 2020)}
\maketitle \thispagestyle{empty}

\noindent \textbf{Abstract}\\
Let $G=(V,E)$ be a connected graph, where $V=\{v_1, v_2, \cdots, v_n\}$.
Let $d_i$ denote the degree of vertex $v_i$.
The ABC matrix of $G$ is defined as $M(G)=(m_{ij})_{n \times n}$,
where $m_{ij}=\sqrt{(d_i + d_j -2)/(d_i d_j)}$
if $v_i v_j \in E$, and 0 otherwise.
The ABC spectral radius of $G$ is the largest eigenvalue of $M(G)$.
In the present paper, we establish two graph perturbations
with respect to ABC spectral radius.
By applying these perturbations, the trees with the third, fourth,
and fifth largest ABC spectral radii are determined.

\baselineskip=0.30in


\baselineskip=0.30in

\section{Introduction}
Let $G=(V,E)$ be a simple connected graph, where $V=\{v_1, v_2, \cdots, v_n\}$.
Let $d_i= d(v_i)$ denote the degree of vertex $v_i$.
$\Delta= \Delta(G) = \max_{1\leq i \leq n} d_i$ is the maximum degree of $G$.
As usual, $S_n$, $P_n$, $C_n$, and $K_n$ denote the star, the path, the cycle,
and the complete graph of order $n$, respectively.   $\mathcal{G}(m,n)$ denotes the set of connected graphs with $n$ vertices and $m$ edges.
$\mathcal{T}_n$ will denote the set of trees of order $n$,
and $\mathcal{T}_n ^{ (\Delta)}= \{ T \in \mathcal{T}_n| \Delta(T) = \Delta \}$.

The atom-bond connectivity (ABC) index of $G$ is defined \cite{b1} as
$ABC(G) = \Sigma_{v_i v_j \in E} f(d_i, d_j)$, where $f(x,y)=\sqrt{ (x+y-2)/(xy) }$.
This topological index turns out to be closely correlated with the heat of formation of alkanes,
and became a hot topic in the last few years (see \cite{b2,b3,b4} and the references therein).
In 2017, Estrada \cite{b5} introduced the ABC matrix of $G$
as $M(G)=(m_{ij})_{n \times n}$,
where $m_{ij}= f(d_i, d_j)$ if $v_i v_j \in E$, and $0$ otherwise.
The chemical background of this matrix was explicated in \cite{b5}.

The eigenvalues of $M = M(G)$ are called the ABC eigenvalues of $G$.
In particular, $\rho_{ABC}(G) = \rho(M)$ is called the ABC spectral radius of $G$,
where $\rho(M)$ is the spectral radius of $M$.
Since $M$ is non-negative, symmetric, and irreducible,
$\rho_{ABC} (G)$ is positive and simple,
and there exists a unique vector $x >0$
such that $\rho_{ABC}(G) = \max_{\|y\|=1} y^T My = x^T Mx$,
which is known as the Perron vector of $M$.

Estrada \cite{b5} first observed that
$\frac{2}{n} ABC(G) \leq \rho_{ABC} (G) \leq \max_{1 \leq i \leq n} M_i$,
with both equalities iff $M_i = \sum_{1 \leq j \leq n} m_{ij}$ is the same, $i=1,2,\cdots,n$.
Recently, Chen \cite{b6} presented another lower bound as
$\rho_{ABC}(G) \geq \sqrt{2(n-R_{-1}(G))/n}$,
where $R_{-1}(G) = \sum_{v_i v_j \in E} \frac{1}{d_i d_j}$.

Chen \cite{b6} further proposed the problem of characterizing graphs
with extremal ABC spectral radius for a given graph class.
Soon, this problem for trees, unicyclic graphs, and connected graphs,
were solved by Chen \cite{b7}, Li et al. \cite{b8}, and Ghorbani et al. \cite{b9}, respectively.

Let $S_{a,b}$ ($a \geq b$ and $a+b = n-2$) denote the double star of order $n$,
and $T_1, T_2, \cdots, T_{10}$ be the trees shown in Figure 1.
\begin{figure}[ht]
  \centering
  \includegraphics[width=5 in]{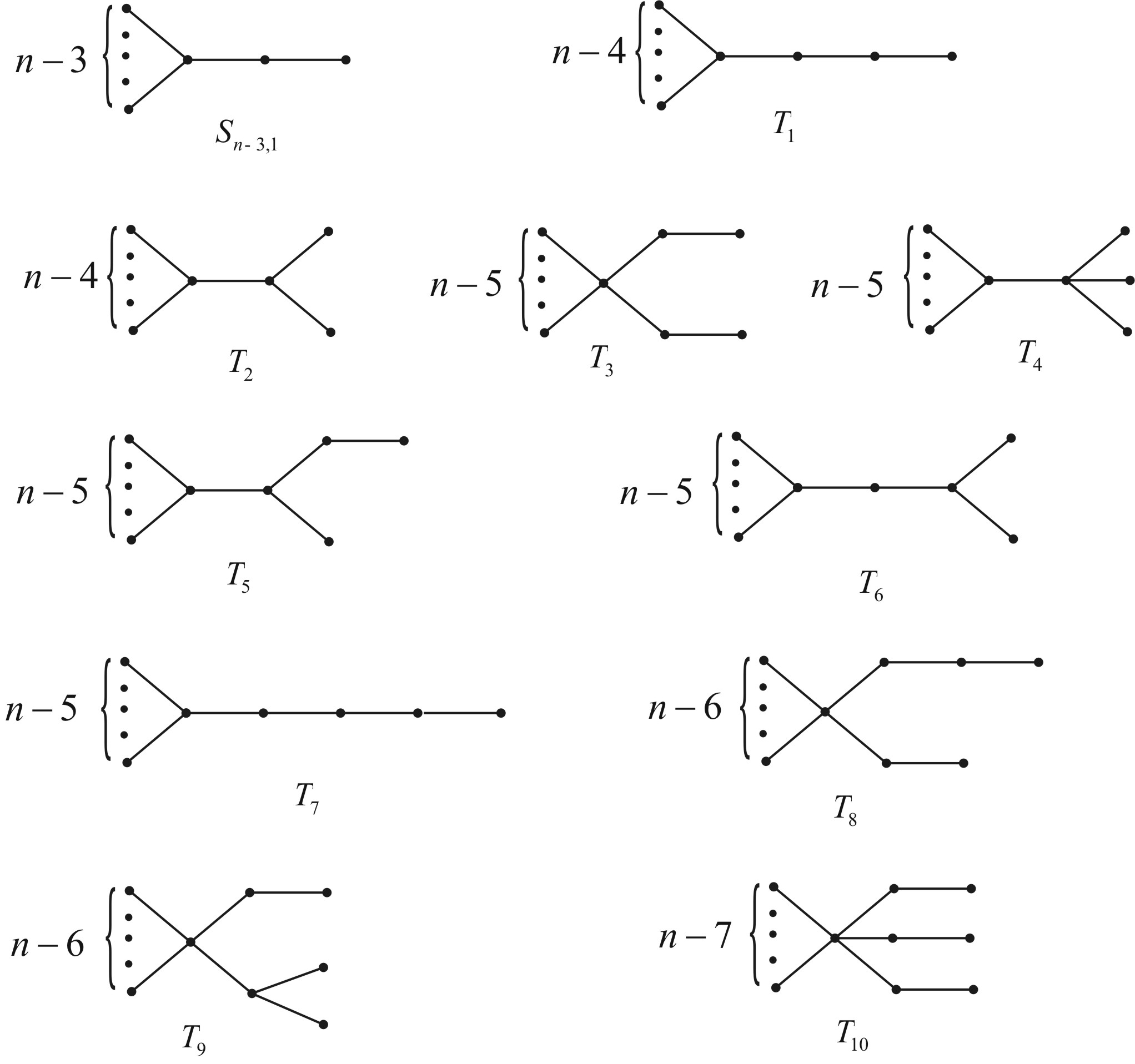}\\
  \caption{The trees $S_{n-3,1}$ and $T_i$, $1 \leq i \leq 10$.}
\end{figure}
Very recently, Lin et al. \cite{b10} established an attainable upper bound of ABC spectral radius as follows.

\noindent \textbf{Lemma 1.1} \cite{b10}\textbf{.} Let $G$ be a connected graph
with $n$ vertices and $m$ edges. Then
$$\rho_{ABC}(G) \leq \sqrt{\Delta(G)+(2m-n+1) / \Delta(G) -2}.$$

Lemma 1.1 reveals that, $\rho_{ABC}(G)$ is large only if $\Delta(G)$ is large.
By applying Lemma 1.1, the unique tree with the second largest ABC spectral radius was determined easily.

\noindent \textbf{Lemma 1.2} \cite{b10}\textbf{.} If $n \geq 4$
and $T \in \mathcal{T}_n - \{S_n, S_{n-3,1} \}$, then
$$\rho_{ABC}(T) < \rho_{ABC}(S_{n-3,1}) < \rho_{ABC}(S_n).$$

Let $\lambda_1(G)$ denote the (adjacent) spectral radius of $G$.
Recall that, while ordering trees in $\mathcal{T}_n$ by their spectral radii,
Lin and Guo \cite{b11} proved the following result.

\noindent \textbf{Theorem 1.3} \cite{b11}\textbf{.} Let $T^{(\Delta)}$
be a tree in $\mathcal{T}_n ^{ (\Delta)}$ and $n \geq 4$. Then
$$\lambda_1 \left( T^{(n-1)} \right) > \lambda_1 \left(T^{(n-2)} \right)
> \cdots \lambda_1 \left( T^{(\lceil \frac{2n}{3} \rceil )} \right)
> \lambda_1 \left(S_{\lceil \frac{2n}{3} \rceil -2, \lfloor \frac{n}{3} \rfloor} \right)
\geq \lambda_1 \left(T^{(k)} \right),$$
where $2 \leq k \leq \lceil \frac{2n}{3} \rceil -1$,
with the equality iff $T^{(k)} \cong S_{\lceil \frac{2n}{3} \rceil -2, \lfloor \frac{n}{3} \rfloor}$.

Naturally, the following question was proposed in \cite{b11}.

\noindent \textbf{Question 1.4.} Let $G_1$ and $G_2$ be two graphs in a subset of $\mathcal{G}(m,n)$.
Is there some integer $l(m,n)$ (depending on $n$ and/or $m$),
such that if $\Delta(G_1) > \Delta(G_2) \geq l(m,n)$,
then $\rho_{ABC}(G_1) > \rho_{ABC}(G_2)$?

This question may be hard to answer at the present.
However, we can attempt to order trees by their ABC spectral radii.
The ordering results may indicate some research directions towards the answer of the question.
In the present paper, we firstly investigate graph perturbations with respect to ABC spectral radius,
and two non-trivial results are obtained.
By applying the perturbations, we determine the trees
with the third, fourth, and fifth largest ABC spectral radii.

\section{Two graph perturbations}
For a graph $G=(V, E)$ with $\{u,v\} \in V$, $G-u$ denotes the subgraph of $G$ induced by $V- \{u\}$.
If $uv \in E$, then $G-uv$ will denote the graph obtained from $G$ by deleting edge $uv$.
Otherwise, if $uv \notin E$, then $G+uv$ will denote the graph obtained from $G$ by adding edge $uv$.

Let $G$ and $H$ be two disjoint graphs, $u \in V(G)$ and $v \in V(H)$.
$G \bigcup H$  will denote the disjoint union of $G$ and $H$.
Let $G(u,v)H = G \bigcup H + uv$. That is, $G(u,v)H$ is the graph obtained
from $G \bigcup H$ by adding edge $uv$.
In particular, if $G$ is non-trivial, $v_0 \in V(G)$, and $P_k = v_1 v_2 \cdots v_k$
and $P_l = v_{-1} v_{-2} \cdots v_{-l}$ are two paths, $k \geq l \geq 0$,
then $G(v_0,v_1)P_k (v_0, v_{-1})P_l$
will be denoted by $G_{k,l}$. If $w \in V(G)$ and $S_{k+1}$ ($S_{l+1}$)
is a star with center $u$ (resp. $v$),
then $G(w,u)S_{k+1}(w,v)S_{l+1}$ will be denoted by $G_{k,l}^{1}$.
$G_{k,l}$ and $G_{k,l}^{1}$ are shown in Figure 2.
\begin{figure}[ht]
  \centering
  \includegraphics[width=4 in]{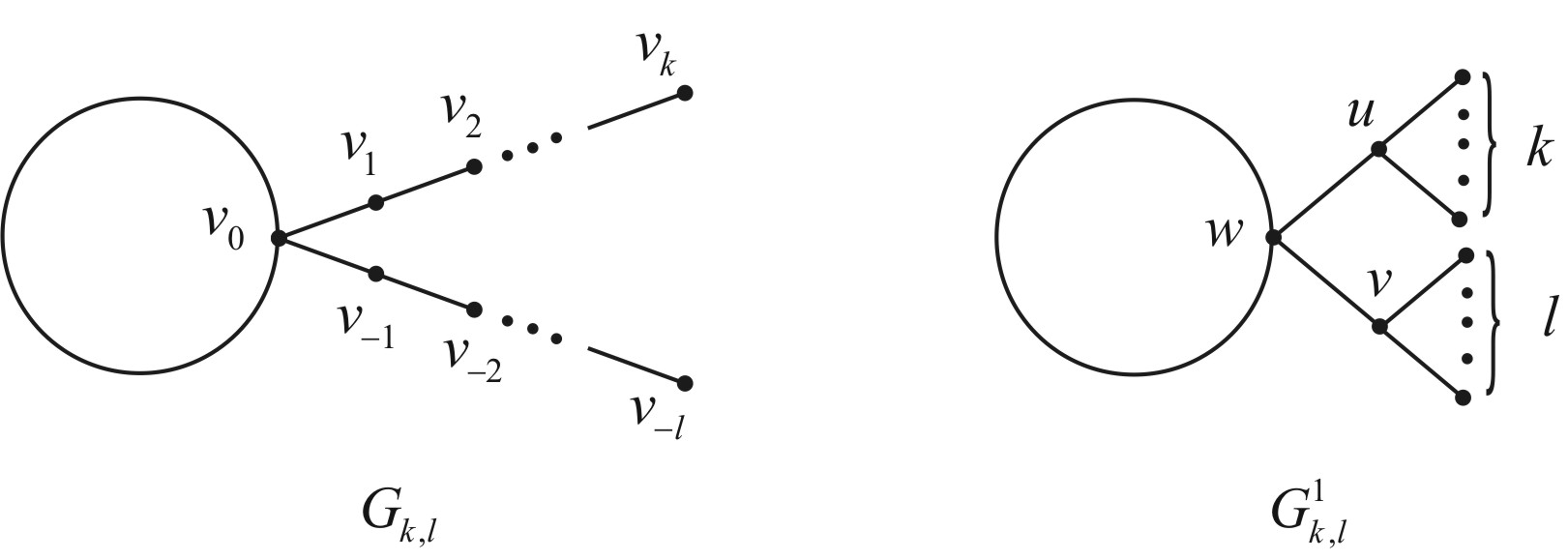}\\
  \caption{The graphs $G_{k,l}$ and $G_{k,l}^1$.}
\end{figure}

Let $P(M,\lambda) = |\lambda I -M|$ be the characteristic polynomial of
the ABC matrix $M = M(G)$ of a graph $G$.
It is both consistent and convenient to define $M(K_1) = \begin{pmatrix} 0 \end{pmatrix}$,
and thus $P \left(M(K_1),\lambda \right) = \lambda$.
Also define $P(\emptyset, \lambda) = 1$,
where $\emptyset$ stands for the virtual square matrix of order 0.
If $H$ is an induced proper subgraph of $G$, then $M_H$ will denote the submatrix of $M$,
which consists of the rows and columns corresponding to the vertices in $V(H)$.
Note that, $M_H$ may be different from $M(H)$, the ABC matrix of graph $H$.
Thus $M_H$ is non-negative and symmetric.
Obviously, it holds that $0 \leq \rho(M_H) < \rho(M) = \rho$ and $P(M_H, \rho) > 0$.

The following result is obvious.

\noindent \textbf{Lemma 2.1.} Let $M = M(G \bigcup H)$. Then $P(M, \lambda) = P(M_G, \lambda) P(M_H, \lambda)$.

\noindent \textbf{Lemma 2.2.} Let $M = M\left( G(u,v)H \right)$. Then
$$P(M, \lambda) = P(M_G, \lambda) P(M_H, \lambda) - f^2 \left( d(u), d(v) \right)
P(M_{G-u}, \lambda) P(M_{H-v}, \lambda).$$

\noindent \textbf{Proof.} Suppose the orders of $G$ and $G(u,v)H$ are $k$ and $n$, respectively.
Label the vertices of $G(u,v)H$ such that $V(G) = \{u = v_1, v_2, \cdots, v_k \}$
and $V(H) = \{v = v_{k+1}, v_{k+2}, \cdots, v_n \}$.

Let $a = f \left( d(u), d(v) \right)$, and

\begin{equation*}
  B = \begin{pmatrix}
    -a &0 &\cdots &0 \\
     0 &0 &\cdots &0 \\
     \cdots &\cdots &\cdots &\cdots \\
     0 &0 &\cdots &0
    \end{pmatrix}_{k \times (n-k)}.
\end{equation*}
Then
\begin{equation*}
  P(M, \lambda) = \begin{vmatrix}
      \lambda I - M_G & B \\
      B^T   &\lambda I - M_H
    \end{vmatrix}.
\end{equation*}

Let $D(j_1, j_2, \cdots, j_k)$ denote the minor of order $k$ of $P(M, \lambda)$,
which consists of the first $k$ rows and the $j_1$-th, $j_2$-th, $\cdots$, $j_k$-th columns,
$1 \leq j_1 < j_2 < \cdots < j_k \leq n$.
Let $C(j_1, j_2, \cdots, j_k)$ be the cofactor of $D(j_1, j_2, \cdots, j_k)$.
It is easily seen that
$$D(1,2,\cdots,k) C(1,2,\cdots,k) = | \lambda I -M_G | |\lambda I -M_H | = P(M_G, \lambda) P(M_H, \lambda),$$
\begin{align*}
  D(2,3,\cdots,k+1) C(2,3,\cdots,k+1) &= (-1)^{k+1} (-a) | \lambda I -M_{G-u} |
                                               \cdot (-a) |\lambda I -M_{H-v} | \\
                                      &= (-1)^{k+1} a^2 P(M_{G-u}, \lambda) P(M_{H-v}, \lambda),
\end{align*}
and $D(j_1, j_2, \cdots, j_k) C(j_1, j_2, \cdots, j_k) = 0$
if $\{j_1, j_2, \cdots, j_k \}$ is not $\{1,2, \cdots,k \}$ or $\{2,3, \cdots, k+1 \}$.
Hence from the Laplace's expansion theorem, we have
\begin{align*}
  P(M, \lambda)&= D(1,2,\cdots,k) C(1,2,\cdots,k) + (-1)^k D(2,3,\cdots,k+1) C(2,3,\cdots,k+1)\\
               &= P(M_G, \lambda) P(M_H, \lambda) - f^2 \left( d(u), d(v) \right)
                                                      P(M_{G-u}, \lambda) P(M_{H-v}, \lambda).
\end{align*}

The proof is thus completed. $\blacksquare$

Without risk of confusion, we denote by $\rho(Q)$ the largest (real) root of a univariate function $
Q(\lambda)$.
With Lemmas 2.1 and 2.2, we are able to give the perturbation of $G_{k,l}^1$ as follows.

\noindent \textbf{Theorem 2.3.} If $k \geq l \geq 1$,
then $\rho_{ABC} \left( G_{k+1,l-1}^1 \right) > \rho_{ABC} \left( G_{k,l}^1 \right)$.

\noindent \textbf{Proof.} By applying Lemmas 2.1 and 2.2,
we compute the characteristic polynomial of $M = M \left( G_{k,l}^1 \right)$ as follows,
where $d = d(w) \geq 2$.
\begin{align*}
  P(M, \lambda)&=\lambda^{k+l} \left[\lambda - \frac{k^2}{(k+1) \lambda} \right]
                               \left[\lambda - \frac{l^2}{(l+1) \lambda} \right] P(M_G, \lambda) -\lambda^{k+l} \\
     &~~~\left\{ \frac{k+d-1}{(k+1)d} \left[ \lambda - \frac{l^2}{(l+1) \lambda}\right]
                             + \frac{l+d-1}{(l+1)d} \left[ \lambda - \frac{k^2}{(k+1) \lambda} \right]
                       \right\} P(M_{G-w}, \lambda) \\
     &= \lambda^{k+l} \left[ \lambda - \frac{k^2}{(k+1) \lambda} \right]
                      \left[ \lambda - \frac{l^2}{(l+1)\lambda} \right]\\
     &~~~\left\{ P(M_G, \lambda) - \frac{\lambda}{d}
                             \left[ \frac{k+d-1}{(k+1) \lambda^2 - k^2} +
                                    \frac{l+d-1}{(l+1) \lambda^2 - l^2} \right]
                                      P(M_{G-w}, \lambda) \right\}.
\end{align*}

Let $\rho = \rho_{ABC}\left( G_{k,l}^1 \right)$.
If $\lambda = \rho$ satisfies $\lambda - k^2/[(k+1)\lambda]$ (that is, $\rho = k/ \sqrt{k+1}$),
then $\rho_{ABC}\left( G_{k+1,l-1}^1 \right) > \rho$ holds obviously.
Otherwise, $\rho > k/ \sqrt{k+1} \geq l/ \sqrt{l+1}$ is the largest (real) root of
$$Q_{k,l}(\lambda) = P(M_G, \lambda) - \frac{\lambda}{d} \left[ \frac{k+d-1}{(k+1) \lambda^2-k^2} + \frac{l+d-1}{(l+1)\lambda^2-l^2} \right] P(M_{G-w}, \lambda).$$

Note that, $P(M_G, \rho) >0$ and $P(M_{G-w}, \rho) >0$.
Thus to prove $\rho_{ABC}\left( G_{k+1,l-1}^1 \right) > \rho $,
it suffices to show $\rho\left( Q_{k+1,l-1}\right) > \rho\left( Q_{k,l} \right) = \rho$.
We will complete the proof by confirming $Q_{k+1,l-1}(\rho) < 0$.

From $Q_{k,l}(\rho) = 0$ we have
$$P(M_G, \rho) = \frac{\rho}{d} \left[ \frac{k+d-1}{(k+1)\rho^2 - k^2} + \frac{l+d-1}{(l+1)\rho^2 - l^2}
  P(M_{G-w}, \rho)  \right].$$
Hence
\begin{align*}
  Q_{k+1,l-1}(\rho) &= \frac{\rho}{d} \left[ \frac{k+d-1}{(k+1)\rho^2 - k^2}
                     + \frac{l+d-1}{(l+1)\rho^2 - l^2} \right] P(M_{G-w}, \rho) \\
                    &~~~ - \frac{\rho}{d} \left[ \frac{k+d}{(k+2)\rho^2 - (k+1)^2}
                     + \frac{l+d-2}{l\rho^2 - (l-1)^2} \right] P(M_{G-w}, \rho),
\end{align*}
\noindent and
\begin{align*}
  &~~~~Q_{k+1,l-1}(\rho) <0\\
  &\Leftrightarrow \frac{k+d-1}{(k+1)\rho^2 - k^2} + \frac{l+d-1}{(l+1)\rho^2 - l^2}
                   <\frac{k+d}{(k+2)\rho^2 - (k+1)^2} + \frac{l+d-2}{l\rho^2 - (l-1)^2}\\
  &\Leftrightarrow \frac{k+d}{(k+2)\rho^2 - (k+1)^2} - \frac{k+d-1}{(k+1)\rho^2 - k^2}
                   >\frac{l+d-1}{(l+1)\rho^2 - l^2} - \frac{l+d-2}{l\rho^2 - (l-1)^2}.
\end{align*}

Let $h(k) = (k+d-1)/[(k+1)\rho^2 - k^2]$, where $k \geq 0$, $d \geq 2$, and $\rho > k/ \sqrt{k+1}$.
Consequently, we will show $h'(k)$ strictly increases with $k > 0$.
By elementary calculations, we have
\begin{align*}
  h''(k) =& \frac{2}{\left[ (k+1)\rho^2 - k^2 \right]^3} \left\{ (k+d-1)(\rho^2 - 2k)^2 \right.\\
          &+ (k+d-1) \left[ (k+1)\rho^2 -k^2\right]
          \left. -(\rho^2 -2k)\left[ (k+1)\rho^2 -k^2 \right] \right\}.
\end{align*}

If $\rho^2-2k \leq 0$, we immediately have $h''(k) > 0$, since $k+d-1 \geq k+1 > 1$
and $(k+1)\rho^2 - k^2 >0$ from $k \geq 1$, $d \geq 2$, and $\rho > k/ \sqrt{k+1}$.
Hence assume $\rho^2 - 2k > 0$.

Since $k \geq 1$, we have $(k+1)\rho^2 -k^2 \geq (k+1) (\rho^2 -2k)$. Hence
\begin{align*}
  h''(k) > 0 &\Leftarrow (k+1)(\rho^2 - 2k)^2 + (k+1)^2(\rho^2 - 2k)
               > (\rho^2 - 2k)\left[(k+1)\rho^2 - k^2 \right]\\
             &\Leftrightarrow (k+1)(\rho^2-2k)+(k+1)^2 > (k+1)\rho^2 -k^2\\
             &\Leftrightarrow 1 > 0.
\end{align*}
By the Lagrange's mean value theorem, there exist two constants $\xi_{k+1}$ and $\xi_l$:
$l-1 < \xi_l < l \leq k < \xi_{k+1} < k+1$, such that
$h(k+1)-h(k) = h'(\xi_{k+1}) > h'(\xi_l) = h(l) - h(l-1)$.
Hence $Q_{k+1,l-1}(\rho) < 0$, and $\rho(Q{k+1,l-1}) > \rho(Q{k,l})$,
which completes the proof. $\blacksquare$

In the rest of this section, we consider the perturbation of $G_{k,l}$
with respect to ABC spectral radius. The relation between $G_{k+1,l-1}$ and $G_{k,l}$,
$k \geq l \geq 1$, is not so good as desired.
For example, the tree $T_2$ in Figure 1 can be regarded as $G_{1,1}$, and $T_1$ as $G_{2,0}$.
We have $\rho_{ABC}(T_2) > \rho_{ABC}(T_1)$ if $n= 6,7,8$,
while $\rho_{ABC}(T_2) < \rho_{ABC}(T_1)$ if $n \geq 9$ (see Lemma 3.4).
According to our numerical experiment, in many cases it holds that
\begin{align*}
\rho_{ABC}(G_{k+l-1,1}) >
  \rho_{ABC}\left( G_{\lceil \frac{k+l}{2} \rceil, \lfloor \frac{k+l}{2} \rfloor} \right)
  & > \rho_{ABC}\left( G_{\lceil \frac{k+l}{2} \rceil -1, \lfloor \frac{k+l}{2} \rfloor +1} \right) \\
  & > \cdots > \rho_{ABC}(G_{k+l-2,2}) > \rho_{ABC}(G_{k+l,0}).
\end{align*}
However, at the present we can have only an almost complete result.

\noindent \textbf{Theorem 2.4.} (1) If $k \geq l \geq 3$,
then $\rho_{ABC}(G_{k,l}) > \rho_{ABC}(G_{k+1,l-1})$.

(2) If $k \geq l \geq 2$, then $\rho_{ABC}(G_{k+l-1,1}) > \rho_{ABC}(G_{k,l})$.

(3) If $k \geq 1$ and all neighbors of $v_0$ in $G$ are of degree 1 or 2,
then $\rho_{ABC}(G_{k,1}) > \rho_{ABC}(G_{k+1,0})$.

To prove Theorem 2.4, we need the following auxiliary result.

\noindent \textbf{Lemma 2.5.} Let $v_0 v_1 \cdots v_k$ be a pendent path
of a connected graph $G$ with $d(v_0) \geq 3$ and $k \geq 2$.
Let $\rho = \rho_{ABC}(G)$ and $\gamma = \left( \rho + \sqrt{\rho^2 -2} \right) / \sqrt{2}$.
Let $x$ be the Perron vector of $M(G)$ with component $x_i$ indexed by $v_i$, $0 \leq i \leq k$.\
Then $x_i = x_0(k+1-i)/(k+1)$ if $\rho = \sqrt{2}$,
and $x_i = x_0 \gamma^i \left( \gamma^{2k+2-2i} -1 \right) / \left( \gamma^{2k+2} -1 \right)$
if $\rho > \sqrt{2}$.
Moreover, $x_i$ strictly increases with $k$, and $x_i < x_0 / \gamma^i$.

\noindent \textbf{Proof.} Since $d(v_0) \geq 3$, $G$ is not a path.
From the proof of the Theorem 2.4 in \cite{b9} we have $\rho \geq \sqrt{2}$,
and $\gamma \geq 1$ is real.

From $\rho x = M(G)x$ we have $\sqrt{2}\rho x_k = x_{k-1}$, and $\sqrt{2}\rho x_i = x_{i-1} + x_{i+1}$
for $1 \leq i \leq k-1$. Extend the sequence $\{ x_i \}_{i=0}^k$ to $\{ x_i \}_{i \geq 0}$
such that the following recurrence equation holds for $i \geq 1$.
\begin{equation}
  \sqrt{2}\rho x_i = x_{i-1} + x_{i+1}.
\end{equation}
Note that, $x_{k+1} = 0$ from $\sqrt{2}\rho x_k = x_{k-1} = x_{k-1} + x_{k+1}$.

It is easily seen that $\gamma$ and $1/ \gamma$ are the roots of the characteristic equation of (2.1).
From the theory of linear recurrence equations,
we will determine the general formula of $x_i$ in the following.

If $\rho =\sqrt{2}$, that is, $\gamma = 1/ \gamma = 1$, we easily get
$$x_i = x_0 \frac{k+1-i}{k+1} < x_0 = \frac{x_0}{\gamma^i},$$
and $x_i$ strictly increases with $k$.

Otherwise, if $\rho > \sqrt{2}$, then $\gamma > 1$,
and there exist constants $a$ and $b$ such that
$x_i = a \gamma^i + b/\gamma^i$.
From the boundary conditions $x_0 = a + b$
and $x_{k+1} = a \gamma^{k+1} + b / \gamma^{k+1} = 0$,
it follows that $a = x_0 / (1-\gamma^{2k+2})$
and $b = x_0 \gamma^{2k+2} / (\gamma^{2k+2} -1)$.
By elementary calculations, we have
$x_i = x_0 \gamma^i \left( \gamma^{2k+2-2i} -1 \right) / \left( \gamma^{2k+2} -1 \right)$.
Thus $x_i < x_0 / \gamma^i$ holds immediately from
$$\frac{\gamma^{2k+2-2i}-1}{\gamma^{2k+2}-1}
  < \frac{\gamma^{2k+2-2i}}{\gamma^{2k+2}} = \frac{1}{\gamma^{2i}}.$$
Finally, we have
$$\frac{\gamma^{2k-2i}-1}{\gamma^{2k}-1} =
\frac{\gamma^{2k+2-2i}-\gamma^2}{\gamma^{2k+2}-\gamma^2}
  < \frac{\gamma^{2k+2-2i}-1}{\gamma^{2k+2}-1}.$$
That is, $x_i$ strictly increases with $k$. $\blacksquare$

\noindent \textbf{Proof of Theorem 2.4.}
(1) Let $M_1$ and $M_2$ be the ABC matrices of $G_{k,l}$ and $G_{k+1,l-1}$, respectively.
Let $y$ be the Perron vector of $M_2$ with component $y_i$ indexed by vertex $v_i$, $1-l \leq i \leq k+1$.
Let $x$ be a vector, whose components are same with $y$, except $x_i$,
which is indexed by $v_i$ of $G_{k,l}$, $-l \leq i \leq k$ and $i \neq 0$.
Set $x_i$, $-l \leq i \leq k$, such that the sequences
$\{x_0, x_1, \cdots, x_k \}$ and $\{x_0, x_{-1}, \cdots, x_{-l}\}$
both satisfy the recurrence equation (2.1).

For convenience, the row of $M_i$ corresponding to vertex $v$ will be denoted by $M_i(v)$, $i=1,2$.
Also denote by $x_v$ and $y_v$ the components of $x$ and $y$ indexed by $v$, respectively.
Based on the Perron-Frobenius theory,
we will complete the proof by confirming $M_1 x \geq \rho x$
and $M_1 x \neq \rho x$.

If $v \neq v_i$, $-l \leq i \leq k$, then $M_1(v)x = M_1(v)y = \rho y_v = \rho x_v$.
It is also easily seen that $M_1(v_i)x = \rho x_i$ for $i \neq 0$.
Hence it remains to confirm $M_1(v_0)x > \rho x_0$.

Obviously,
$\sqrt{2} [M_1(v_0)x - \rho x_0] = \sqrt{2} [M_1(v_0)x - \rho y_0 ] = x_1 + x_{-1} - y_1 -y_{-1}$.

\noindent \textbf{Case 1.} $\rho = \sqrt{2}$. Since $k \geq l \geq 3$,
we have $l(l+1) < (k+1)(k+2)$. Hence from Lemma 2.5 we have
\begin{align*}
  [x_1 + x_{-1} - y_1 - y_{-1}]/x_0 &= \frac{k}{k+1} + \frac{l}{l+1} - \frac{k+1}{k+2} - \frac{l-1}{l} \\
                                    &= \left[\frac{l}{l+1} - \frac{l-1}{l}\right] - \left[\frac{k+1}{k+2}
                                                                                      -\frac{k}{k+1}\right]\\
                                    &= \frac{l}{l(l+1)} - \frac{1}{(k+1)(k+2)} \\
                                    &> 0.
\end{align*}

\noindent \textbf{Case 2.} $\rho > \sqrt{2}$.
Then $\gamma = \left( \rho + \sqrt{\rho^2 -2} \right) / \sqrt{2} >1$.
From Lemma 2.5 we have
\begin{align*}
  [x_1 + x_{-1} - y_1 -y_{-1}]/x_0 & = \gamma\left(\frac{\gamma^{2l}-1}{\gamma^{2l+2}-1}
                                         -\frac{\gamma^{2l-2}-1}{\gamma^{2l}-1} \right)
                                      -\gamma\left(\frac{\gamma^{2k+2}-1}{\gamma^{2k+4}-1}
                                         - \frac{\gamma^{2k}-1}{\gamma^{2k+2}-1} \right) \\
   &= \gamma (\gamma^2-1)^2 \left[\frac{\gamma^{2l-2}} {(\gamma^{2l+2}-1)(\gamma^{2l}-1)} -
                                    \frac{\gamma^{2k}}{ (\gamma^{2k+4}-1)(\gamma^{2k+2}-1)} \right]
\end{align*}
It is easily seen that
$$h(k) = \frac{\gamma^{2k}}{ (\gamma^{2k+4}-1)(\gamma^{2k+2}-1)}
       = \left[\left(\gamma^{k+4}- 1/ \gamma^k \right)\left(\gamma^{k+2}-1/ \gamma^k \right) \right]^{-1}$$
strictly decreases with $k>0$. Hence $ x_1 + x_{-1} - y_1 -y_{-1} >0$, which completes the proof.

(2) Let $M_1$ and $M_2$ be the ABC matrices of $G_{k+l-1,1}$ and $G_{k,l}$, respectively,
and $\rho = \rho_{ABC}(G_{k,l})$. Let $d \geq 1$ be the degree of $v_0$ in $G$.
We distinguish the following two cases.

\noindent \textbf{Case 1.} $\rho = \sqrt{2}$. By properly label the vertices of $G_{k+l-1,1}$,
from the basic properties of $f(x,y)$ (for example, see \cite{b3}),
$M_1$ has a submatrix whose entries are all not less than those
of the following non-negative and symmetry matrix:
\begin{equation*}
  \bar{M} = \begin{pmatrix}
     0              &\sqrt{\frac{d+1}{d+2}} &\sqrt{\frac{1}{2}} &\sqrt{\frac{1}{d+2}} &\cdots &\sqrt{\frac{1}{d+2}}\\
\sqrt{\frac{d+1}{d+2}}        &0                 &0                     &0            &\cdots &0 \\
\sqrt{\frac{1}{2}}  &0        &0                 &0                     &0            &\cdots &0 \\
\sqrt{\frac{1}{d+2}}&0        &0                 &0                     &0            &\cdots &0 \\
     \cdots                   &\cdots            &\cdots                &\cdots       &\cdots &\cdots\\
\sqrt{\frac{1}{d+2}}&0        &0                 &0                     &0            &\cdots &0
   \end{pmatrix}_{(d+3)\times(d+3)}.
\end{equation*}
It is easily seen that
$$\rho(M_1) > \rho(\bar{M}) = \sqrt{\frac{d+1}{d+2} +\frac{1}{2} +\frac{d}{d+2}}
= \sqrt{\frac{5d+4}{2d+4}}> \sqrt{2} =\rho.$$

\noindent \textbf{Case 2.} $\rho > \sqrt{2}$.
Let $y$ be the Perron vector of $M_2$ with component $y_i$ indexed by $v_i$ of $G_{k,l}$,
$-l \leq i \leq k$. Let $x>0$ be a vector, whose components are same with $y$,
except $x_1, x_2, \cdots, x_{k+l-1}$ and $x_{-1}$,
which are indexed by $v_1, v_2, \cdots, v_{k+l-1}$ and $v_{-1}$ of $G_{k+l-1,1}$, respectively.
Set $x_i$, $-1 \leq i \leq k+l-1$, such that the sequence $\{x_0, x_1, \cdots, x_{k+l-1}\}$
satisfies the recurrence equation (2.1) and $x_{-1} = f(1,d+2)x_0 / \rho$.
Note that $x_0 = y_0$.

It is easily seen that $M_1(v)x = \rho x_v$ if $v \neq v_0$.
Finally, from Lemma 2.5 we have
\begin{align*}
  M_1(v_0)x - \rho x_0 &= M_1(v_0)x - \rho y_0 \\
    &= \sqrt{\frac{1}{2}} x_1 + f(1,d+2) x_{-1} -\sqrt{\frac{1}{2}} y_1 - \sqrt{\frac{1}{2}} y_{-1} \\
    &> \sqrt{\frac{d+1}{d+2}}\frac{x_0}{\rho} - \frac{x_0}{\sqrt{2}\gamma}.
\end{align*}
Hence
\begin{align*}
  M_1(v_0)x > \rho x_0 &\Leftarrow \frac{\sqrt{d+1}}{\sqrt{d+2}\rho} \geq \frac{1}{\sqrt{2} \gamma} \\
    &\Leftrightarrow 2(d+1)\gamma^2 \geq (d+2)\rho^2 \\
    &\Leftrightarrow (d+1)\left(\rho+\sqrt{\rho^2-2} \right)^2 \geq (d+2)\rho^2\\
    &\Leftrightarrow 2(d+1)\rho^2 + (d+1)\left(2\rho\sqrt{\rho^2-2}-2\right) \geq (d+2)\rho^2\\
    &\Leftarrow d\rho^2 >0.
\end{align*}

(3) Let $M_1 = M(G_{k,1})$ and $M_2 = M(G_{k+1,0})$.
Suppose $u_1, u_2, \cdots, u_d$ are the neighbors of $v_0$ in $G$.
Let $x$ be the Perron vector of $M_2$ with component $x_i$ indexed by $v_i$,
and $y_j$ by $u_j$ in $G_{k+1,0}$, $0 \leq i \leq k+1$ and $1 \leq j \leq d$.
Regard $G_{k,1}$ as $G_{k+1,0} - v_k v_{k+1} + v_0 v_{k+1}$.
From the basic properties of $f(x,y)$ we have
\begin{align*}
  x^T (M_1-M_2)x &= f(1,d+2)x_0 x_{k+1} -\sqrt{\frac{1}{2}}x_k x_{k+1} \\
    &~~~+ \sum_{i=1}^{d} \left[f(d(u_i), d+2) - f(d(u_i), d+1) \right]x_0 y_i \\
    &\geq f(1,d+2)x_0 x_{k+1} -\sqrt{\frac{1}{2}}x_k x_{k+1}\\
    &>0.
\end{align*}

The proof is thus completed. $\blacksquare$

\section{Ordering trees by their ABC spectral radii}
By applying Theorems 2.3 and 2.4, we will prove the following result in this section.

\noindent \textbf{Theorem 3.1.} If $n \geq 10$ and $T \in \mathcal{T}_n - \{S_n, S_{n-3,1}, T_1, T_2, T_3\}$, then $$\rho_{ABC}(T) < \rho_{ABC}(T_3) < \rho_{ABC}(T_2) < \rho_{ABC}(T_1)< \rho_{ABC}(S_{n-3,1}) <\rho_{ABC}(S_n).$$

Bearing Lemma 1.2 in mind, the proof of Theorem 3.1 will be completed by the following six lemmas.

\noindent \textbf{Lemma 3.2.} If $n \geq 10$, $\Delta \leq n-5$, and $T \in \mathcal{T}_n^{(\Delta)}$,
then $\rho_{ABC}(T) < \sqrt{n-5} < \rho_{ABC}(T_3)$.

\noindent \textbf{Proof.} Let
\begin{equation*}
  \bar{M} =\begin{pmatrix}
      0 &\sqrt{\frac{n-4}{n-3}} &\sqrt{\frac{n-4}{n-3}} &\cdots & &\sqrt{\frac{n-4}{n-3}} \\
      \sqrt{\frac{n-4}{n-3}} &0 &0 &\cdots &0 \\
      \sqrt{\frac{n-4}{n-3}} &0 &0 &\cdots &0 \\
      &\cdots &\cdots &\cdots &\cdots &\cdots \\
      \sqrt{\frac{n-4}{n-3}} &0 &0 &\cdots &0 \\
    \end{pmatrix}_{(n-3) \times (n-3)}.
\end{equation*}

It is easily seen that $\bar{M}$ is a submatrix of $M(T_3)$. Hence
$$\rho_{ABC}(T_3) \geq \rho(\bar{M}) =\sqrt{\frac{(n-4)^2}{n-3} } > \sqrt{n-5}.$$

On the other hand, from Lemma 1.1 and $n \geq 10$ we have
$$\rho_{ABC}(T) \leq \sqrt{n-5+(n-1)/(n-5) -2} < \sqrt{n-5}.~\blacksquare$$

Thus to prove Theorem 3.1, it suffices to consider the trees in
$\mathcal{T}_n^{n-3} = \{T_1, T_2, T_3\}$ and $\mathcal{T}_n^{n-4} = \{T_i| 4 \leq i \leq 10 \}$.

\noindent \textbf{Lemma 3.3.} If $n \geq 6$,
then $\rho_{ABC}(T_2) > \rho_{ABC}(T_3)$.

\noindent \textbf{Proof.} Immediately from Theorem 2.3. $\blacksquare$

\noindent \textbf{Lemma 3.4.} $\rho_{ABC}(T_1) < \rho_{ABC}(T_2)$ if $6 \leq n \leq 8$,
and $\rho_{ABC}(T_1) > \rho_{ABC}(T_2)$ if $n \geq 9$.

\noindent \textbf{Proof.} For $6 \leq n \leq 12$, the conclusion can be confirmed easily.
Hence assume $n \geq 13$. From Lemmas 2.1 and 2.2, we easily get
$$P(M(T_1), \lambda) = \lambda^{n-3} \left[\lambda-\frac{(n-4)^2}{(n-3)\lambda} \right](\lambda^2-1)
                         -\frac{1}{2}\lambda^{n-4} \left(\lambda^2 - \frac{1}{2} \right),$$
$$P(M(T_2), \lambda) = \lambda^{n-2} \left[\lambda-\frac{(n-4)^2}{(n-3)\lambda} \right]
                  \left(\lambda-\frac{4}{3 \lambda}\right)-\frac{n-2}{3(n-3)}\lambda^{n-2}.$$

Let
$$Q_1(\lambda) = \lambda \left[\lambda-\frac{(n-4)^2}{(n-3)\lambda} \right](\lambda^2-1)
                         -\frac{1}{2} \left(\lambda^2 - \frac{1}{2} \right),$$

$$Q_2(\lambda) = \left[\lambda-\frac{(n-4)^2}{(n-3)\lambda} \right]
                          \left(\lambda-\frac{4}{3 \lambda}\right)-\frac{n-2}{3(n-3)}.$$

Obviously, $\rho = \rho_{ABC}(T_2) = \rho(Q_2) > \sqrt{5}$.
From $Q_2(\rho) = 0$ we have
$$\rho-\frac{(n-4)^2}{(n-3)\rho} = \frac{n-2}{3(n-3)\left(\rho- \frac{4}{3\rho}\right)},$$
and
\begin{align*}
 Q_1(\rho) &= \rho(\rho^2-1) \frac{n-2}{3(n-3)\left(\rho- \frac{4}{3\rho}\right)}
                - \frac{1}{2} \left(\rho^2 - \frac{1}{2} \right) \\
           &< \frac{(n-2)\rho(\rho^2-1)}{3(n-3)\left(\rho- \frac{4}{3\rho}\right)}
                - \frac{1}{2} (\rho^2-1).
\end{align*}

Since $n \geq 13$, we have
\begin{align*}
 Q_1(\rho) < 0 &\Leftarrow \frac{(n-2)\rho}{3(n-3)\left(\rho- \frac{4}{3\rho}\right)} \leq \frac{1}{2} \\
          &\Leftarrow \frac{11\rho}{30\left(\rho- \frac{4}{3\rho}\right)} \leq \frac{1}{2} \\
          &\Leftrightarrow \rho \geq \sqrt{5},
\end{align*}
and it follows that $\rho_{ABC}(T_1) =\rho(Q_1) > \rho = \rho_{ABC}(T_2)$. $\blacksquare$

\noindent \textbf{Lemma 3.5.} If $n \geq 10$, then $\rho_{ABC}(T_{10}) < \rho_{ABC}(T_9) < \rho_{ABC}(T_4) <\sqrt{n-5} < \rho_{ABC}(T_3)$.

\noindent \textbf{Proof.} From Theorem 2.3 we immediately have
$\rho_{ABC}(T_{10}) < \rho_{ABC}(T_9) < \rho_{ABC}(T_4)$.
It remains to prove $\rho_{ABC}(T_4) < \sqrt{n-5}$, since $\rho_{ABC}(T_3) > \sqrt{n-5}$ from Lemma 3.2.
From Lemmas 2.1 and 2.2, we easily get
$$P(M(T_4),\lambda) = \lambda^{n-2} \left[\lambda- \frac{(n-5)^2}{(n-4)\lambda}\right]
                                    \left(\lambda- \frac{9}{4\lambda}\right) -\lambda^{n-2}\frac{n-2}{4(n-4)}.$$

Let
$$Q(\lambda) = \left[\lambda- \frac{(n-5)^2}{(n-4)\lambda}\right]
                 \left(\lambda- \frac{9}{4\lambda}\right) - \frac{n-2}{4(n-4)}.$$
It is easily seen that $Q(\lambda)$ strictly increases with $\lambda \geq \sqrt{n-5}$, and
\begin{align*}
  Q(\sqrt{n-5}) >0 &\Leftrightarrow \left[\sqrt{n-5}- \frac{(n-5)^2}{(n-4)\sqrt{n-5}}\right]
                 \left(\sqrt{n-5}- \frac{9}{4\sqrt{n-5}}\right) > \frac{n-2}{4(n-4)} \\
    &\Leftrightarrow \left[n-5-\frac{(n-5)^2}{n-4}\right] \left(n-5-\frac{9}{4}\right) > \frac{(n-2)(n-5)}{4(n-4)} \\
    &\Leftrightarrow \left[(n-4)(n-5)-(n-5)^2\right] \left(n-5-\frac{9}{4}\right) > \frac{(n-2)(n-5)}{4} \\
    &\Leftrightarrow n >9.
\end{align*}
Hence $\rho = \rho(Q) < \sqrt{n-5}$, which completes the proof. $\blacksquare$

\noindent \textbf{Lemma 3.6.} If $n \geq 8$, then $\rho_{ABC}(T_8) < \rho_{ABC}(T_7) < \rho_{ABC}(T_6)
< \sqrt{n-5} < \rho_{ABC}(T_3)$.

\noindent \textbf{Proof.} From Theorem 2.4, $\rho_{ABC}(T_8) < \rho_{ABC}(T_7) < \rho_{ABC}(T_6)$ holds immediately.
It remains to show $\rho_{ABC}(T_6) < \sqrt{n-5}$ from Lemma 3.2.

From Lemmas 2.1 and 2.2 we easily get
$$P(M(T_6),\lambda) = \lambda^{n-2}\left[\lambda-\frac{(n-5)^2}{(n-4)\lambda}\right]
                      \left(\lambda-\frac{11}{6\lambda}\right)
                        -\lambda^{n-3}\left(\frac{\lambda}{2} -\frac{4}{6\lambda}\right).$$

Let
$$Q(\lambda) = \lambda \left[\lambda-\frac{(n-5)^2}{(n-4)\lambda}\right]
                       \left(\lambda-\frac{11}{6\lambda}\right) -\left(\frac{\lambda}{2} -\frac{4}{6\lambda}\right).$$

It is easily seen that $\rho = \rho_{ABC}(T_6) > \sqrt{n-6} \geq \sqrt{3} > \sqrt{17/6}$.
Hence we have
$$\frac{\rho}{2} -\frac{4}{6\rho} < \frac{3}{4}\left(\rho- \frac{11}{6\rho}\right),$$
and
\begin{align*}
  0 = Q(\rho) &= \rho \left[\rho-\frac{(n-5)^2}{(n-4)\rho}\right]
                       \left(\rho-\frac{11}{6\rho}\right) -\left(\frac{\rho}{2} -\frac{4}{6\rho}\right) \\
              &> \rho \left[\rho-\frac{(n-5)^2}{(n-4)\rho}\right]
                       \left(\rho-\frac{11}{6\rho}\right) -\frac{3}{4}\left(\rho-\frac{11}{6\rho}\right) \\
              &= \left(\rho-\frac{11}{6\rho}\right)\left[\rho^2 -\frac{(n-5)^2}{n-4} -\frac{3}{4}\right].
\end{align*}

Obviously, $\rho-11/(6\rho) > 0$, hence $\rho^2 -(n-5)^2/(n-4) -3/4 <0$.
Since $n \geq 8$, it follows that
$$\rho < \sqrt{\frac{4n^2-37n+88}{4(n-4)}} \leq \sqrt{\frac{4n^2-36n+80}{4(n-4)}} = \sqrt{n-5}.~\blacksquare$$

\noindent \textbf{Lemma 3.7.} If $n \geq 7$, then $\rho_{ABC}(T_5) < \sqrt{n-5} < \rho_{ABC}(T_3)$.

\noindent \textbf{Proof.} For $n=7,8$, the conclusion can be confirmed easily..
Hence assume $n \geq 9$. From Lemmas 2.1 and 2.2 we easily get
$$P(M(T_5),\lambda) = \lambda^{n-5}\left[\lambda -\frac{(n-5)^2}{(n-4)\lambda}\right]
                      \left(\lambda^4 -\frac{5}{3}\lambda^2 +\frac{1}{3}\right)
                       -\frac{n-3}{3(n-4)}\lambda^{n-4} \left(\lambda^2 -\frac{1}{2} \right).$$

Let
$$Q(\lambda) = \left[\lambda -\frac{(n-5)^2}{(n-4)\lambda}\right]
                      \left(\lambda^4 -\frac{5}{3}\lambda^2 +\frac{1}{3}\right)
                       -\frac{n-3}{3(n-4)} \lambda \left(\lambda^2 -\frac{1}{2}\right),$$
and $\rho = \rho_{ABC}(T_5) = \rho(Q) > \sqrt{n-6} \geq \sqrt{3} > \sqrt{\frac{17}{6}}$.
We have
$$\rho^4-\frac{5}{3}\rho^2+\frac{1}{3} > \rho^2\left(\rho^2-\frac{5}{3}\right) > \frac{1}{2}\rho^2 \left(\rho^2-\frac{1}{2}\right).$$

On the other hand, since $n \geq 9$ we have
$$\frac{n-3}{3(n-4)} \rho \left(\rho^2 -\frac{1}{2}\right) \leq \frac{2}{5} \rho \left(\rho^2 -\frac{1}{2}\right).$$

Hence
\begin{align*}
 0 = Q(\rho) &= \left[\rho -\frac{(n-5)^2}{(n-4)\rho}\right]
                      \left(\rho^4 -\frac{5}{3}\rho^2 +\frac{1}{3}\right)
                       -\frac{n-3}{3(n-4)} \rho \left(\rho^2 -\frac{1}{2} \right) \\
             &> \left[\rho -\frac{(n-5)^2}{(n-4)\rho}\right] \frac{1}{2}\rho^2 \left(\rho^2-\frac{1}{2}\right)
                      - \frac{2}{5} \rho \left(\rho^2 -\frac{1}{2}\right) \\
             &= \frac{1}{2}\rho \left(\rho^2-\frac{1}{2}\right)
                   \left[\rho^2 -\frac{(n-5)^2}{n-4} - \frac{4}{5}\right],
\end{align*}
and we have $\rho^2 -(n-5)^2/(n-4) - 4/5 <0$. Therefore
$$\rho < \sqrt{\frac{(n-5)^2}{n-4} + \frac{4}{5}} = \sqrt{\frac{5n^2-46n+109}{5(n-4)}}
      \leq \sqrt{\frac{5n^2-45n+100}{5(n-4)}} = \sqrt{n-5}.~\blacksquare$$

\section{Further discussions}
In the present paper, we give perturbations with respect to ABC spectral radius
for $G_{k,l}$ and $G_{k,l}^1$. By applying these perturbations,
we determine the trees of order $n \geq 10$ with the third, fourth, and fifth
largest ABC spectral radii.
Though it is possible to extend the ordering,
we leave it a task in the future,
especially after giving the complete perturbation of $G_{k,l}^1$.

From Theorem 3.1 we know that, for two trees $T_1$ and $T_2$ of order $n \geq 10$,
if $\Delta(T_1) > \Delta(T_2) \geq n-4$, then $\rho_{ABC}(T_1) > \rho_{ABC}(T_2)$.
However, comparing with Theorem 1.3, which concerns the spectral radius of trees,
this result is still too trivial.
It seems necessary to characterize the extremal trees in $\mathcal{T}_n^{(\Delta)}$,
at least for those whose maximum degree is large.
Based on the two graph perturbations, for $\Delta \geq \lceil \frac{n}{2} \rceil$,
we guess in $\mathcal{T}_n^{(\Delta)}$, the double star $S_{\Delta-1, n-\Delta-1}$
or the bloom $B_{n,n-\Delta-1}$ maximizes the ABC spectral radius,
while $S(n,n-\Delta-1, 2\Delta-n+1)$ minimizes the ABC spectral radius.
$B_{n,n-\Delta-1}$ and $S(n,n-\Delta-1, 2\Delta-n+1)$ are shown in Figure 3.

\begin{figure}[ht]
  \centering
  \includegraphics[width=4 in]{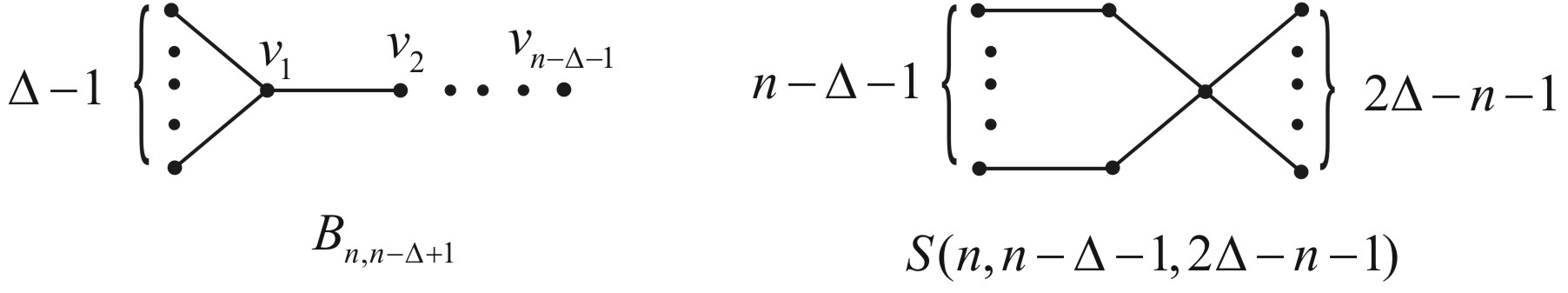}\\
  \caption{The trees $B_{n,n-\Delta-1}$ and $S(n,n-\Delta-1, 2\Delta-n+1)$, $\Delta \geq \lceil \frac{n}{2} \rceil$.}
\end{figure}

Finally, for the graph $G$ concerned in Lemma 2.5, we guess $\rho_{ABC}(G) > \sqrt{2}$.
It may be interesting to characterize connected graphs with small ABC spectral radius,
because after all, till now the known lower bounds are somehow trivial.
Hence we end this paper with the following conjecture.

\noindent \textbf{Conjecture 4.1.} If $G$ is a connected graph of order
$n \geq 4$ and $\rho_{ABC}(G) \leq \sqrt{2}$, then
$G \in \{P_n, C_n, S_4\}$.

\end{document}